\documentclass[12pt,a4paper]{amsart}
\usepackage{amssymb,amscd}

\theoremstyle{plain}
\newtheorem{theorem}{Theorem}[section]
\newtheorem{lemma}[theorem]{Lemma}

\newtheorem{corollary}[theorem]{Corollary}

\theoremstyle{definition}

\newtheorem{remarks}[theorem]{Remarks}
\newtheorem{example}[theorem]{Example}

\numberwithin{equation}{section}

\newcommand\bk{{\Bbbk}}

\newcommand\bA{{\mathbb A}}

\newcommand\bG{{\mathbb G}}

\newcommand\bP{{\mathbb P}}

\newcommand\cA{{\mathcal A}}

\newcommand\aff{\operatorname{aff}}

\newcommand\diag{\operatorname{diag}}

\newcommand\Aut{\operatorname{Aut}}

\newcommand\End{\operatorname{End}}

\newcommand\Ker{\operatorname{Ker}}

\newcommand\Stab{\operatorname{Stab}}

\title{The structure of normal algebraic monoids}

\author{Michel Brion and Alvaro Rittatore}  
\address{Universit\'e de Grenoble I\\
D\'epartement de Math\'ematiques\\
Institut Fourier, UMR 5582 du CNRS\\
38402 Saint-Martin d'H\`eres Cedex, France}
\email{Michel.Brion@ujf-grenoble.fr}

\address{Facultad de Ciencias\\ 
Universidad de la Rep\'ublica\\
Igu\'a 4225\\
11400 Montevideo, Uruguay}
\email{alvaro@cmat.edu.uy}

\begin{document}
 
\begin{abstract}
We show that any normal algebraic monoid is an extension of an abelian
variety by a normal affine algebraic monoid. This extends (and builds
on) Chevalley's structure theorem for algebraic groups.
\end{abstract}

\maketitle

\section{Introduction}
\label{sec:introduction}

A classical theorem of Chevalley asserts that any connected algebraic
group is an extension of an abelian variety by a connected affine
algebraic group. In this note, we obtain an analogous result for
normal algebraic monoids. This reduces their structure to that of more 
familiar objects: abelian varieties, and affine (equivalently, linear)
algebraic monoids. The latter have been extensively investigated, see
the expositions \cite{Pu88, Re05}. 

To state Chevalley's theorem and our analogue in a precise way, we
introduce some notation. We consider algebraic varieties and algebraic
groups over an algebraically closed field $\bk$ of arbitrary
characteristic. By a variety, we mean a separated integral scheme of
finite type $X$ over $\bk$; by a point of $X$, we mean a closed
point. An algebraic group is a smooth group scheme of finite type over
$\bk$.

Let $G$ be a connected algebraic group, then there exists a unique
connected normal affine algebraic subgroup $G_{\aff}$ such that the quotient
$G/G_{\aff}$ is an abelian variety. In other words, we have an exact
sequence of connected algebraic groups
\begin{equation}\label{eqn:extension}
\CD 
1 @>>> G_{\aff} @>>> G  @>{\alpha_G}>> \cA(G) @>>> 0
\endCD
\end{equation}
where $G_{\aff}$ is affine and $\cA(G)$ is projective (since the
group $\cA(G)$ is commutative, its law will be denoted additively).
It follows that the morphism $\alpha_G$ is affine; hence the variety
$G$ is quasi-projective (see \cite{Co02} for these developments and
for a modern proof of Chevalley's theorem).

Next, let $M$ be an irreducible \emph{algebraic monoid}, i.e., an
algebraic variety over $\bk$ equipped with a morphism $M \times M \to M$
(the \emph{product}, denoted simply by $(x,y) \mapsto xy$) which is
associative and admits an identity element $1$. Denote by $G = G(M)$
the group of invertible elements of $M$. The \emph{unit group} $G$ is
known to be a connected algebraic group, open in $M$ (see
\cite[Thm.~1]{Ri98}). 

Let $G_{\aff} \subseteq G$ be the associated affine group, and
$M_{\aff}$ the closure of $G_{\aff}$ in $M$. Clearly, $M_{\aff}$ is an
irreducible algebraic monoid with unit group $G_{\aff}$. By
\cite[Thm.~2]{Ri06}, it follows that $M_{\aff}$ is affine. Also, note
that
\begin{equation}
\label{eqn:GM}
M = G M_{\aff} = M_{\aff} G
\end{equation}
as follows from the completeness of $G/G_{\aff} = \cA(G)$
(see Lemma \ref{lem:map} for details). 

We may now state our main result, which answers a question raised by
D.~A.~Timashev (see the comments after Thm.~17.3 in \cite{Ti06}):

\begin{theorem}\label{thm:main}
Let $M$ be an irreducible algebraic monoid with unit group $G$. If the
variety $M$ is normal, then $\alpha_G : G \to \cA(G)$ 
extends to a morphism of algebraic monoids $\alpha_M : M \to \cA(G)$.
Moreover, the morphism $\alpha_M$ is affine, and its scheme-theoretic
fibers are normal varieties; the fiber at $1$ equals $M_{\aff}$.
\end{theorem}

In loose words, any normal algebraic monoid is an extension of an
abelian variety by a normal affine algebraic monoid. 

For nonsingular monoids, Theorem \ref{thm:main} follows
immediately from Weil's extension theorem: any rational map from a
nonsingular variety to an abelian variety is a morphism. However, this
general result no longer holds for singular varieties. Also, the
normality assumption in Theorem \ref{thm:main} cannot be omitted, as
shown by Example \ref{ex:nonnormal}.

Some developments and applications of the above theorem are
presented in Section 2. The next section gathers a number of
auxiliary results to be used in the proof of that theorem, given in
Section 4. The final Section 5 contains further applications of our
structure theorem to the classification of normal algebraic monoids,
and to their faithful representations as endomorphisms of
homogeneous vector bundles on abelian varieties.

\medskip

\noindent {\bf Acknowledgements.} The second author would like to
thank the Institut Fourier for its hospitality; his research was also
partially supported by the Fondo Clemente Estable, Uruguay (FCE-10018).

\section{Some applications}
\label{sec:applications}

With the notation and assumptions of Theorem \ref{thm:main}, 
observe that $\alpha_M$ is equivariant with respect to the action of
the group $G \times G$ on $M$ via 
$$
(g_1,g_2) \cdot m = g_1 m g_2^{-1},
$$ 
and its action on $\cA(G)$ via 
$$
(g_1,g_2) \cdot a = \alpha(g_1) - \alpha(g_2) + a.
$$
Since the latter action is transitive, $\alpha_M$ is a 
$G \times G$-homogeneous fibration. In particular, all fibers are
isomorphic, and $\alpha_M$ is faithfully flat. 

Also, each irreducible component of the closed subset 
$M_{\aff} \setminus G_{\aff} \subset M_{\aff}$ is of codimension $1$,
since $G_{\aff}$ is affine. Together with the above observation, it
follows that the same holds for the set $M \setminus G$ of non-units
in $M$:

\begin{corollary}\label{cor:noninvertible}
Each irreducible component of $M \setminus G$ has codimension $1$ in $M$.
\qed
\end{corollary}

Next, we obtain an intrinsic characterization of the morphism
$\alpha_M$. To state it, recall from \cite{Se58} that any variety $X$
admits an \emph{Albanese morphism}, i.e., a universal morphism to an
abelian variety.

\begin{corollary}
\label{cor:albanese}
$\alpha_M$ is the Albanese morphism of the variety $M$.
\end{corollary}

\begin{proof}
Let $f : M \to A$ be a morphism (of varieties) to an abelian
variety. Composing $f$ with a translation of $A$, we may assume that
$f(1) = 0$. Then the restriction $f\vert_G : G \to A$ is a morphism of
algebraic groups by \cite[Lem.~2.2]{Co02}. So $f(G_{\aff})$ equals
$0$ by \cite[Lem.~2.3]{Co02}. It follows that 
$f\vert_G = \varphi \circ \alpha_G$, where $\varphi : \cA(G) \to A$ is
a morphism of algebraic groups. Hence $f$ equals 
$\varphi \circ \alpha_M$, since both morphisms have the same
restriction to the open subset $G$.
\end{proof}

Also, since the morphism $\alpha_M : M \to \cA(G)$ is affine and the
variety $\cA(G)$ is projective (see e.g. \cite[Thm.~7.1]{Mi86}), 
we obtain the following:

\begin{corollary}
\label{cor:quasi-projective}
$M$ is quasi-projective.
\qed
\end{corollary}

Another application of Theorem \ref{thm:main} concerns the set
$$
E(M) = \{ e\in M ~\vert~ e^2 = e\}
$$
of \emph{idempotents}. Indeed, since $\alpha_M$ is a morphism of
monoids and the unique idempotent of $\cA(G)$ is the origin, we
obtain:

\begin{corollary}
\label{cor:idempotents}
$E(M) = E(M_{\aff})$. 
\qed
\end{corollary}

Next, recall that a monoid $N$ is said to be \emph{regular}\/ if 
given any $x \in N$, there exists $y \in N$ such that $x = x y x$.

\begin{corollary}\label{cor:regular}
$M$ is regular if and only if $M_{\aff}$ is regular.
\end{corollary}

\begin{proof}
If $M_{\aff}$ is regular, then so is $M$ by (\ref{eqn:GM}).
Conversely, assume that $M$ is regular. Let $x \in M_{\aff}$ and write
$x = x y x$, where $y \in M$. Then we obtain: $\alpha_M(y) = 0$, so
that $y \in M_{\aff}$.
\end{proof}

By \cite[Thm.~13]{Pu82}, every regular irreducible affine algebraic
monoid $N$ is \emph{unit regular}, i.e., given any $x \in N$, there
exists $y \in G(N)$ such that $x = x y x$; equivalently, 
$N = G(N) E(N)$. Together with Corollary \ref{cor:regular}, this
implies:

\begin{corollary}\label{cor:unit}
If $M$ is regular, then it is unit regular.
\qed
\end{corollary} 

Finally, we show that Theorem \ref{thm:main} does not extend to
arbitrary irreducible algebraic monoids:

\begin{example}\label{ex:nonnormal}
Let $A$ be an abelian variety. Then 
$$
M := A \times \bA^1
$$ 
is a commutative nonsingular algebraic monoid via the product
$$
(x,y) \; (x',y') = (x + x', yy'),
$$
with unit group 
$$
G:= A \times \bG_m
$$ 
and kernel $A \times \{ 0 \}$. The morphism $\alpha_G$ is the first
projection $A \times \bG_m \to A$; likewise, the Albanese morphism of
the variety $M$ is just the first projection
$$
p: A \times \bA^1 \to A.
$$

Next, let $F \subset A$ be a non-trivial finite subgroup. Let $M'$ be
the topological space obtained from $M$ by replacing the closed subset
$A \times \{ 0 \}$ with the quotient $A/F \times \{ 0 \}$; in other
words, each point $(x + f,0)$  (where $x \in A$ and $f \in F$) is
identified with the point $(x, 0)$. Denote by
$$
q : M \to M'
$$
the natural map. We claim that $M'$ has a structure of a
irreducible, non-normal, commutative algebraic monoid with unit group 
$G$, such that $q$ is a morphism of monoids; furthermore, 
$\alpha_G: G \to A$ does not extend to a morphism $M' \to A$.

Indeed, one readily checks that $M'$ carries a unique product such that 
$q$ is a morphism of monoids. Moreover, the restriction $q\vert_G$ is
an isomorphism onto $G(M')$, and we have a commutative square
$$
\CD
M \ @>{q}>> M' \\
@V{ p }VV  @V{\alpha}VV \\
A @>{ }>>  A/F \\
\endCD
$$
where $\alpha$ is equivariant with respect to the action of $A$ on $M'$
via the product of $M'$, and the natural action of $A$ on $A/F$. Let 
$$
N := \alpha^{-1}(0),
$$
then the set $N$ is the image under $q$ of the subset
$F \times \bA^1 \subset A \times \bA^1$.
So $N$ is a union of copies of the affine line, indexed by the finite
set $F$, and glued along the origin. Hence $N$ is a reduced affine
scheme, and its product (induced by the product of $M'$) is a
morphism: $N$ is a connected, reducible affine algebraic
monoid. Furthermore, the natural map
$$
A \times^F N \to M'
$$
is clearly an isomorphism of monoids, and the left-hand side is also
an algebraic monoid. This yields the desired structure of algebraic
monoid on $M'$. The map $q$ is induced from the natural map
$F \times \bA^1  \to N$, which is a morphism; hence so
is $q$. Finally, the projection $p: A \times \bG_m \to A$ cannot
extend to a morphism $M' \to A$: such a morphism would be
$A$-equivariant, and hence restrict to an $A$-equivariant morphism
$A/F \times \{ 0 \} \cong A/F \to A$, which is impossible.
This completes the proof of the claim.

Alternatively, this claim follows from a general result concerning
the existence of pinched schemes (see \cite[Thm.~5.4]{Fe03}). Indeed, 
$M'$ is obtained by pinching the quasi-projective variety $M$ along
its closed subset $A \times \{ 0 \} \cong A$ via the finite morphism 
$A \to A/F$, in the terminology of \cite{Fe03}.

One easily checks that $\alpha$ is the Albanese morphism of the
variety $M'$, and $q$ is its normalization. Moreover,
$$
M'_{\aff} = q( \{ 0 \} \times \bA^1) \cong \bA^1
$$ 
is strictly contained in $N$, and is nonsingular whereas $M'$ is 
non-normal. 

Note finally that $M'$ is weakly normal, i.e., any finite bijective
birational map from a variety to $M'$ is an isomorphism. Thus, Theorem
\ref{thm:main} does not extend to weakly normal monoids.
\end{example}

\section{Auxiliary results}
\label{sec:auxiliary}

We consider a connected algebraic group $G$ and denote by $Z^0$
its connected center regarded as a closed reduced subscheme of $G$,
and hence as a connected algebraic subgroup. 

\begin{lemma}\label{lem:product}
{\rm (i)} The scheme-theoretic intersection $Z^0 \cap G_{\aff}$
contains $Z^0_{\aff}$ as a normal subgroup, and the quotient
$(Z^0 \cap G_{\aff})/Z^0_{\aff}$ is a finite group scheme.

\smallskip

\noindent
{\rm (ii)} The product map $Z^0 \times G_{\aff} \to G$ factors through
an isomorphism
$$
(Z^0 \times G_{\aff})/(Z^0 \cap G_{\aff}) \cong G,
$$
where $Z^0 \cap G_{\aff}$ is embedded in $Z^0 \times G_{\aff}$ as a
normal subgroup scheme via the identity map on the first factor and
the inverse map on the second factor.

\smallskip

\noindent
{\rm (iii)} The natural map 
$Z^0/(Z^0 \cap G_{\aff}) \to G/G_{\aff} = \cA(G)$ is an isomorphism of
algebraic groups.
\end{lemma}

This easy result is proved in \cite[Sec.~1.1]{Br06} under the
assumption that $\bk$ has characteristic zero; the general case follows
by similar arguments. We will also need the following result, see
e.g. \cite[Sec.~1.2]{Br06}:

\begin{lemma}\label{lem:action}
Let $G$ act faithfully on an algebraic variety $X$. Then the isotropy
subgroup scheme of any point of $X$ is affine.
\qed
\end{lemma}

Next we consider an irreducible algebraic monoid $M$ with unit group
$G$. If $M$ admits a zero element, then this point is fixed by $G$
acting by left multiplication, and this action is faithful. Thus, $G$
is affine by Lemma \ref{lem:action}. Together with
\cite[Thm.~2]{Ri06}, this yields:

\begin{corollary}\label{cor: fixed}
Any irreducible algebraic monoid having a zero element is affine.
\qed
\end{corollary}

Returning to an arbitrary irreducible algebraic monoid $M$, recall
that an \emph{ideal} of $M$ is a subset $I$ such that 
$M I M \subseteq I$.

\begin{lemma}\label{lem:kernel}
{\rm (i)} $M$ contains a unique closed $G\times G$-orbit, which is
also the unique minimal ideal: the kernel $\Ker(M)$.

\smallskip

\noindent
{\rm (ii)} If $M$ is affine, then $\Ker(M)$ contains an idempotent.
\end{lemma}

\begin{proof}
(i) is part of \cite[Thm.~1]{Ri98}. For (ii), see e.g. \cite[p.~35]{Re05}. 
\end{proof}

\begin{lemma}\label{lem:map}
{\rm (i)} Let $Z^0 \cap G_{\aff}$ act on $Z^0 \times M_{\aff}$ by
multiplication on the first factor, and the inverse map composed with
left multiplication on the second factor. Then the quotient 
$$
Z^0 \times^{Z^0 \cap G_{\aff}} M_{\aff} := 
(Z^0 \times M_{\aff})/(Z^0 \cap G_{\aff})
$$ 
has a unique structure of an irreducible algebraic monoid such that
the quotient map is a morphism of algebraic monoids. Moreover,
$$
G\bigl(Z^0 \times^{Z^0 \cap G_{\aff}} M_{\aff}\bigr) = 
Z^0 \times^{Z^0 \cap G_{\aff}} G_{\aff} \cong G.
$$
Regarded as a $G$-variety via left multiplication, 
$Z^0 \times^{Z^0 \cap G_{\aff}} M_{\aff}$ is naturally isomorphic to
the quotient
$$
G \times^{G_{\aff}} M_{\aff} := (G \times M_{\aff})/G_{\aff},
$$
where the action of $G_{\aff}$ on $G \times M_{\aff}$ is defined as above.

\smallskip

\noindent
{\rm (ii)} The product map $Z^0 \times M_{\aff} \to M$ factors
uniquely through a morphism of algebraic monoids
\begin{equation}\label{eqn:map}
\pi: Z^0 \times^{Z^0 \cap G_{\aff}} M_{\aff} \to M.
\end{equation}
Moreover, $\pi$ is birational and proper.

\smallskip

\noindent
{\rm (iii)} $M = Z^0 M_{\aff}$ and $\Ker(M) = Z^0 \Ker(M_{\aff})$.
\end{lemma}

\begin{proof}
(i) and the first assertion of (ii) are straightforward.
Also, the restriction of $\pi$ to the unit group is an isomorphism,
and hence $\pi$ is birational. To show the properness, observe that
$\pi$ factors as a closed immersion
$$
Z^0 \times^{Z^0 \cap G_{\aff}} M_{\aff} \to 
Z^0 \times^{Z^0 \cap G_{\aff}} M
$$
(induced by the inclusion map $M_{\aff} \to M$), followed by an
isomorphism
$$
Z^0 \times^{Z^0 \cap G_{\aff}} M \to (Z^0/Z^0 \cap G_{\aff}) \times M
\cong \cA(G) \times M
$$
given by $(z,m)\mapsto \bigl(z(Z^0 \cap G_{\aff}), zm\bigr)$, followed
in turn by the projection 
$$
\cA(G) \times M \to M
$$
which is proper, since $\cA(G)$ is projective.

(iii) By (ii), $\pi$ is surjective, i.e., the first equality holds. 
For the second equality, note that $Z^0 \Ker(M_{\aff})$ is closed in
$M$ since $\pi$ is proper, and is a unique orbit of $G \times G$ by
Lemma \ref{lem:product}.
\end{proof}

\section{Proof of the main result}
\label{sec:proof}

We begin by showing the following result of independent interest:

\begin{theorem}\label{thm:map}
The morphism 
$$
\pi: Z^0 \times^{Z^0 \cap G_{\aff}} M_{\aff}  = G \times^{G_{\aff}} M 
\to M
$$
is an isomorphism for any normal irreducible monoid $M$.
\end{theorem}

\begin{proof}
We proceed through a succession of reduction steps.

1) It suffices to show that $\pi$ is finite for any irreducible
algebraic monoid $M$ (possibly non-normal). Indeed, the desired
statement follows from this, in view of Lemma \ref{lem:map}(ii) and
Zariski's Main Theorem.

2) Since $\pi$ is proper, it suffices to show that its fibers
are finite. But the points of $M$ where the fiber of $\pi$ is finite
form an open subset (by semicontinuity), which is stable under the
action of $G \times G$. Thus, it suffices to check the finiteness of
the fiber at some point of the unique closed $G \times G$-orbit,
$\Ker(M)$. By Lemmas \ref{lem:kernel} (ii) and \ref{lem:map} (iii), 
$\Ker(M)$ contains an idempotent $e \in \Ker(M_{\aff})$. So we are
reduced to showing that the set $\pi^{-1}(e)$ is finite.

3) Consider the $Z^0$-orbit $Z^0 e$ and its inverse image under
$\pi$,
$$
\pi^{-1}(Z^0 e) \cong
Z^0 \times^{Z^0 \cap G_{\aff}}(Z^0 e \cap M_{\aff}).
$$
It suffices to check that the map
$$
p : Z^0 \times^{Z^0 \cap G_{\aff}}(Z^0 e \cap M_{\aff}) \to Z^0 e 
$$
(the restriction of $\pi$) has finite fibers. Since $p$ is surjective,
it suffices in turn to show that 
$Z^0 \times^{Z^0 \cap G_{\aff}}(Z^0 e \cap M_{\aff})$ and 
$Z^0 e$ are algebraic groups of the same dimension, and $p$ is a
morphism of algebraic groups.

4) Since $e$ is idempotent and $Z^0$ is a central subgroup of $G$,
the orbit $Z^0 e$ (regarded as a locally closed, reduced
subscheme of $M$) is a connected algebraic group under the product of
$M$, with identity element $e$. Moreover, the intersection 
$Z^0 e \cap M_{\aff}$ (also regarded as a locally closed,
reduced subscheme of $M$) is a closed submonoid of $Z^0 e$,
with the same identity element $e$. By \cite[Exer.~3.5.1.2]{Re05}, 
it follows that $Z^0 e \cap M_{\aff}$ is a subgroup of $Z^0 e$. Hence
$Z^0 \times^{Z^0 \cap G_{\aff}}(Z^0 e \cap M_{\aff})$
is an algebraic group as well, and clearly $p$ is a morphism of
algebraic groups. 

5) It remains to show that
\begin{equation}\label{eqn:dim}
\dim Z^0 \times^{Z^0 \cap G_{\aff}}(Z^0 e \cap M_{\aff})
= \dim (Z^0 e).
\end{equation}

We first analyze the left-hand side. Since $Z^0 e \cap M_{\aff}$
is a quasi-affine algebraic group, it is affine. But the maximal
connected affine subgroup of $Z^0 e$ is $Z^0_{\aff} e$,
since $Z^0 e \cong Z^0/ \Stab_{Z^0}(e)$ as groups. It follows that 
$Z^0_{\aff} e \subseteq Z^0 e \cap M_{\aff}$
is the connected component of the identity. Hence
$$
\dim Z^0 \times^{Z^0 \cap G_{\aff}}(Z^0 e \cap M_{\aff}) = 
\dim (Z^0) - \dim (Z^0 \cap G_{\aff}) 
+ \dim(Z^0_{\aff} e).
$$
But $\dim(Z^0 \cap G_{\aff}) = \dim(Z^0_{\aff})$ by Lemma
\ref{lem:product}, and 
$$
\dim(Z^0_{\aff} e) = \dim (Z^0_{\aff}) - \dim \Stab_{Z^0_{\aff}}(e)
$$
so that
$$
\dim Z^0 \times^{Z^0 \cap G_{\aff}}(Z^0 e \cap M_{\aff}) = 
\dim (Z^0) - \dim \Stab_{Z^0_{\aff}}(e).
$$

On the other hand,
$$
\dim (Z^0 e) = \dim (Z^0) - \dim \Stab_{Z^0}(e),
$$ 
and $\Stab_{Z^0}(e)$ is affine by Lemma \ref{lem:action}. Hence 
$$
\dim \Stab_{Z^0}(e) = \dim \Stab_{Z^0_{\aff}}(e).
$$
This completes the proof of Equation (\ref{eqn:dim}) and, in turn, of
the finiteness of $\pi$.
\end{proof}

We may now prove Theorem \ref{thm:main}. Observe that the projection
$$
M \cong Z^0 \times^{Z^0 \cap G_{\aff}} M_{\aff} \to Z^0 /(Z^0 \cap G_{\aff})
\cong \cA(G)
$$
yields the desired extension $\alpha_M$ of $\alpha_G$. Clearly, the
scheme-theoretic fiber of $\alpha_M$ at $1$ equals $M_{\aff}$.

To show that the morphism $\alpha_M$ is affine, oberve that $M_{\aff}$
is an affine variety equipped with an action of the affine group
scheme $Z^0 \cap G_{\aff}$. Thus, $M_{\aff}$ admits a closed
equivariant immersion into a $(Z^0 \cap G_{\aff})$-module $V$. Then
$Z^0 \times^{Z^0 \cap G_{\aff}} M_{\aff}$ (regarded as a variety over
$\cA(G)$) admits a closed $Z^0$-equivariant immersion into 
$Z^0 \times^{Z^0 \cap G_{\aff}} V$, the total space of a vector bundle
over $\cA(G)$.

Finally, to show that the variety $M_{\aff}$ is normal, consider its
normalization $\widetilde{M_{\aff}}$, an affine variety
where $Z^0 \cap G_{\aff}$ acts such that the normalization map 
$f : \widetilde{M_{\aff}} \to M_{\aff}$ is equivariant. This defines a
morphism
$$
\widetilde{\pi} : Z^0 \times^{Z^0 \cap G_{\aff}} \widetilde{M_{\aff}} \to M
$$
which is still birational and finite. So $\widetilde{\pi}$ is an
isomorphism by Zariski's Main Theorem; it follows that $f$ is an
isomorphism as well.
\qed

\section{Classification and faithful representation}
\label{sec:classification}

We begin by reformulating our main results as a classification theorem 
for normal algebraic monoids:

\begin{theorem}\label{thm:equivalence}
The category of normal algebraic monoids is equivalent to the category
having as objects the pairs $(G,N)$, where $G$ is a connected
algebraic group and $N$ is a normal affine algebraic monoid with unit
group $G_{\aff}$. 

The morphisms from such a pair $(G,N)$ to a pair $(G',N')$ are the
pairs $(\varphi,\psi)$, where $\varphi : G \to G'$ is a morphism of
algebraic groups and $\psi: N \to N'$ is a morphism of algebraic
monoids such that $\varphi \vert_{G_{\aff}} = \psi \vert_{G_{\aff}}$.
\end{theorem}

\begin{proof}
By Theorem \ref{thm:map}, any normal irreducible algebraic $M$ with
unit group $G$ is determined by the pair $(G,M_{\aff})$ up to
isomorphism. Conversely, any pair $(G,N)$ as in the above
statement yields a normal algebraic monoid 
$$
M := G \times^{G_{\aff}} N 
$$
together with isomorphisms $G \to G(M)$ and $N \to M_{\aff}$, as
follows from Lemmas \ref{lem:product} and \ref{lem:map}. 

Next, consider a morphism of normal algebraic monoids 
$$
f: M \to M'.
$$
Clearly, $f$ restricts to a morphism of algebraic groups 
$$
\varphi : G(M) \to G(M').
$$ 
Moreover, the universal property of the Albanese maps $\alpha_M$,
$\alpha_{M'}$ yields a commutative diagram
$$
\CD
M @>{f}>> M' \\
@V{\alpha_M}VV  @V{\alpha_{M'}}VV \\
\cA\bigl(G(M)\bigr) @>{\alpha_f}>> \cA\bigl(G(M')\bigr), \\
\endCD
$$
where $\alpha_f$ is a morphism of varieties such that 
$\alpha_f(0) = 0$, and hence a morphism of abelian varieties (see
e.g. \cite[Cor.~3.6]{Mi86}). In turn, this yields a morphism of
algebraic monoids
$$
\psi : M_{\aff} = \alpha_M^{-1}(0) \to \alpha_{M'}^{-1}(0) = M'_{\aff}
$$
such that $\varphi\vert_{G(M)_{\aff}} = \psi\vert_{G(M)_{\aff}}$.
Conversely, any such pair $(\varphi,\psi)$ defines a morphism $f$, as
follows from Theorem \ref{thm:map} again.
\end{proof}

\begin{remarks} (i) The irreducible affine algebraic monoids having a
prescribed unit group $G$ are exactly the affine equivariant
embeddings of the homogeneous space $(G \times G)/\diag(G)$, by
\cite[Prop.~1]{Ri98}. In the case where $G$ is reductive, such
embeddings admit a combinatorial classification, see \cite{Ri98} and
\cite{Ti06}.

\smallskip

\noindent
(ii) The normality assumption in Theorem \ref{thm:equivalence}
cannot be omitted: with the notation of Example \ref{ex:nonnormal},
the monoids $M$ and $M'$ yield the same pair $(A \times \bG_m,\bA^1)$;
but they are not isomorphic as varieties, since the image of their
Albanese map is $A$, resp.~$A/F$.
\end{remarks}

Next, we obtain faithful representations of normal algebraic monoids
as endomorphisms of vector bundles over abelian varieties. For this,
we need additional notation and some preliminary observations.

Let $A$ be an abelian variety, and 
$$
p : E \to A
$$
a vector bundle. Observe that $p$ is the Albanese morphism of the
variety $E$ (as follows e.g. from \cite[Cor.~3.9]{Mi86}). Thus, any
morphism of varieties $f : E \to E$ fits into a commutative square
$$
\CD
E @>{f}>> E \\
@V{p}VV @V{p}VV \\
A @>{\alpha(f)}>> A, \\
\endCD
$$
where $\alpha(f)$ is a morphism of varieties as well. By
\cite[Cor.~3.9]{Mi86} again, $\alpha(f)$ is the composition of a
translation of $A$ with an endomorphism of the abelian variety $A$. 

We say that $f$ is an \emph{endomorphism} (resp.~an
\emph{automorphism}) of $E$, if $\alpha(f)$ is the translation 
$$
t_a : A \to A, \quad x \mapsto a +x
$$ 
for some $a = a(f)\in A$, and the induced maps on fibers 
$$
f_x: E_x \to E_{a+x} \quad (x \in A)
$$ 
are all linear (resp.~linear isomorphisms). 

Clearly, the endomorphisms of $E$ form a monoid under composition,
denoted by $\End(E)$; its unit group $\Aut(E)$ consists of the
automorphisms. The map
\begin{equation}\label{eqn:alb}
\alpha : \End(E) \to A, \quad f \mapsto a(f)
\end{equation}
is a morphism of monoids, and its fiber at a point $a \in A$ is
isomorphic to the set of morphisms of vector bundles from
$E$ to $t_a^* E$ (a finite-dimensional $k$-vector space). 

In particular, the fiber at $0$ is the monoid $\End_A(E)$ of
endomorphisms of $E$ regarded as a vector bundle over $A$. Moreover,
$\End_A(E)$ is a finite-dimensional $k$-algebra; in particular, an
irreducible affine algebraic monoid. Its unit group $\Aut_A(E)$ is the
kernel of the restriction of $\alpha$ to $\Aut(E)$.

The vector bundle $E$ is called \emph{homogeneous} if the
restriction map $\Aut(E) \to A$ is surjective; equivalently, 
$E \cong t_a^* E$ for any $a \in A$. 

For example, a line bundle is homogeneous if and only if it is
algebraically equivalent to $0$ (see \cite[Sect.~9]{Mi86}). More
generally, the homogeneous vector bundles are exactly the direct sums
of vector bundles of the form $L \otimes F$, where $L$ is an
algebraically trivial line bundle, and $F$ admits a filtration by
sub-vector bundles such that the associated graded bundle is trivial
(see \cite[Thm.~4.17]{Mu78}). 

We are now in a position to state:

\begin{theorem}\label{thm:rep}
{\rm (i)} Let $p: E \to A$ be a homogeneous vector bundle over an
abelian variety. Then $\End(E)$ has a structure of a nonsingular
irreducible algebraic monoid such that its action on $E$ is
algebraic. Moreover, the Albanese morphism of $\End(E)$ is the map
of (\ref{eqn:alb}), so that $\End(E)_{\aff} = \End_A(E)$.

\smallskip

\noindent
{\rm (ii)} Any normal irreducible algebraic monoid $M$ is isomorphic to
a closed submonoid of $\End(E)$, where $E$ is a homogeneous vector
bundle over the Albanese variety of $M$.
\end{theorem}

\begin{proof}
(i) We claim that $\Aut(E)$ has a structure of a connected algebraic
group such that its action on $E$ is algebraic and the map 
$\Aut(E) \to A$ is a (surjective) morphism of algebraic groups. 

Indeed, any $f \in \Aut(E)$ extends uniquely to an
automorphism of the projective completion $\bP(E \oplus O_A)$
(regarded as an algebraic variety), where $O_A$ denotes the trivial
bundle of rank $1$ over $A$. Moreover, every automorphism of 
$\bP(E \oplus O_A)$ induces an automorphism of $A$, the Albanese
variety of $\bP(E \oplus O_A)$. It follows that $\Aut(E)$ may be
identified to the group $G$ of automorphisms of $\bP(E \oplus O_A)$
that induce translations of $A$, and commute with the action of the 
multiplicative group $\bG_m$ by multiplication on fibers of $E$.
Clearly, $G$ is contained in the connected automorphism group 
$\Aut^0 \bP(E \oplus O_A)$ (a connected algebraic group) as a closed
subgroup; hence, $G$ is an algebraic group. Moreover, the exact sequence
$$
1 \to \Aut_A(E) \to \Aut(E)  \to A \to 0
$$
implies that $G$ is connected; this completes the proof of our claim.

This exact sequence also implies that 
$\Aut(E)_{\aff} =  \Aut_A(E)$. Hence the natural map 
$$
\pi: \Aut(E) \times^{\Aut_A(E)} \End_A(E) \to \End(E)
$$
is bijective, since the map $\alpha: \End(E) \to A$ is a
$\Aut(E)$-homogeneous fibration. This yields a structure of algebraic
monoid on $\End(E)$, which clearly satisfies our assertions.

(ii) By \cite[Thm.~3.15]{Pu88}, the associated monoid $M_{\aff}$ is
isomorphic to a closed submonoid of $\End(V)$, where $V$ is a
vector space of finite dimension over $\bk$. In particular, $V$ is a
rational $G_{\aff}$-module. We may thus form the associated vector
bundle
$$
p : E := G \times^{G_{\aff}} V = Z^0 \times^{Z^0 \cap G_{\aff}} V 
\to \cA(G).
$$
Since the action of $\cA(G)$ on itself by translations lifts to the
action of $G$ on $E$, then $E$ is homogeneous. Moreover, one easily
checks that the product action of $Z^0 \times M_{\aff}$ on 
$Z^0 \times V$ yields a faithful action of 
$M = Z^0 \times^{Z^0 \cap G_{\aff}} M_{\aff}$ on $E$.
\end{proof}

\end{document}